# Fixed point theorems without continuity in metric vector spaces

Jinlu Li

Department of Mathematics
Shawnee State University
Portsmouth Ohio 45662
USA

**Abstract**

In this oaper, we prove some fixed point theorems in metric vector spaces, in which the continuity is not required for the considered mappings to satisfy. We provide some concrete examples to demonstrate these theorems. We also give some counterexamples to show that every condition in the theorems is necessary for the considered mapping to have a fixed point.

## 1. Introduction

In fixed point theory, a topological vector space $X$ is said to have the fixed point property if any continuous self-mapping on an arbitrary given nonempty compact and convex subset of $X$ has a fixed point. In the concept of fixed point property, a certain type of continuity of the considered mapping is absolutely required (see [2, 3, 5, 6]). In this paper, we prove a fixed point theorem in metric vector spaces. In this theorem, the continuity is not a necessary condition for the considered mapping to satisfy, which is replaced by a certain type of convexity.

After K. Fan proved the so called Fan-KKM theorem in 1961, it has been widely applied in fixed point theory, optimization theory, variational inequalities and other fields (see [4]). The theorems of this paper are also proved by using Fan-KKM theorem with two different ways. Since Fan-KKM Theorem has more than one different versions, for easy reference, we briefly review the version of Fan-KKM Theorem used in this note (see Fan [1] and Park [4]).

Let $C$ be a nonempty convex subset of a vector space $X$. A set-valued mapping $G: C \to 2^X \setminus \{\emptyset\}$ is called a KKM mapping if, for any finite subset $\{x_1, x_2, \ldots, x_n\}$ of $C$, we have

$$\text{co}\{x_1, x_2, \ldots, x_n\} \subseteq \bigcup_{1 \leq i \leq n} G(x_i),$$

where $\text{co}\{x_1, x_2, \ldots, x_n\}$ denotes the convex hull of $\{x_1, x_2, \ldots, x_n\}$.

**Fan-KKM Theorem**. *Let $C$ be a nonempty closed convex subset of a Hausdorff topological vector space $X$ and let $G: C \to 2^X \setminus \{\emptyset\}$ be a KKM mapping with closed values. If there exists a point $x^* \in C$ such that $G(x^*)$ is a compact subset of $C$, then*

$$\bigcap_{x \in C} G(x) \neq \emptyset.$$

## 2. The first fixed point theorem in this paper

### 2.1 The first fixed point theorem in this paper

**Theorem 1.** *Let $(X, d)$ be a metric vector space and let $C$ be a nonempty closed and convex subset of $X$. Let $f: C \to C$ be a single-valued mapping. Suppose that $f$ satisfies the following conditions*:

(a) $f: C \to C$ is onto;

($b_1$) *For any finite subset $\{x_1, x_2, \ldots, x_n\} \subseteq C$ and for arbitrary convex combination $u$ of $x_1, x_2, \ldots, x_n$ with $u = \sum_{i=1}^{n} \alpha_i x_i$, in which $\alpha_1, \alpha_2, \ldots, \alpha_n$ are positive with $\sum_{i=1}^{n} \alpha_i = 1$, we have*
$$\max\{d(f(x_j), u) - d(x_j, u): j = 1, 2, \ldots, n\} \geq 0;$$

($c_1$) *There is $x^* \in C$ such that the following subset of $C$ is compact*
$$\{y \in C: d(x^*, y) \leq d(f(x^*), y)\}.$$

*Then $f$ has a fixed point*.

*Proof.* We define a set-valued mapping $G: C \to 2^C$ by
$$G(x) = \{y \in C: d(x, y) \leq d(f(x), y)\}, \text{ for every } x \in C.$$

For every $x \in C$, $G(x)$ is nonempty because $x \in G(x)$. $G(x)$ is a closed subset of $C$ since, for every fixed $x \in C$, $d(f(x), \cdot) - d(x, \cdot)$ is a continuous function on $C$ (on $X$).

Next, we show that $G$ is a KKM mapping. To this end, for any finite subset $\{x_1, x_2, \ldots, x_n\} \subseteq C$, let $u$ be a (linear) convex combination of $x_1, x_2, \ldots, x_n$. We can suppose that there are positive numbers $\alpha_1, \alpha_2, \ldots, \alpha_n$ with $\sum_{i=1}^{n} \alpha_i = 1$ such that $u = \sum_{i=1}^{n} \alpha_i x_i$. By condition ($b_1$) in this theorem, we have
$$\max\{d(f(x_j), u) - d(x_j, u): j = 1, 2, \ldots, n\} \geq 0;$$

It implies that there must be an integer $k$ with $1 \leq k \leq n$ such that
$$d(x_k, u) \leq d(f(x_k), u)$$
That is,
$$\sum_{i=1}^{n} \alpha_i x_i = u \in G(x_k) \subseteq \bigcup_{1 \leq j \leq n} G(x_j).$$

It implies that $G: C \to 2^C$ is a KKM mapping with nonempty closed values in $C$. By condition ($c_1$) and by using Fan-KKM Theorem, we obtain $\bigcap_{x \in C} G(x) \neq \emptyset$. Then, taking any $y_0 \in \bigcap_{x \in C} G(x)$, we have

$$d(x, y_0) \leq d(f(x), y_0), \text{ for every } x \in C. \tag{1}$$

By condition (a) in this theorem, for the arbitrarily selected $y_0 \in \cap_{x \in C} G(x) \subseteq C$ satisfying (1), there is $x_0 \in C$ such that $f(x_0) = y_0$. Substituting $x$ by $x_0$ in (1) gets

$$d(x_0, f(x_0)) = d(x_0, y_0) \leq d(f(x_0), y_0) = d(y_0, y_0) = 0.$$

It implies that $x_0$ is a fixed point of $f$, which proves this theorem. □

**Corollary 3.** *Let $(X, d)$ be a metric vector space and let $C$ be a nonempty compact and convex subset of $X$. Let $f: C \to C$ be a single-valued mapping. If $f$ satisfies conditions* (a) *and* (b$_1$) *in Theorem 1, then $f$ has a fixed point.*

### 2.2 Examples regarding to Theorem 1

In the following examples, we always take $(X, d) = (\mathbb{R}, |\cdot|)$ and let $C$ be a closed interval of $\mathbb{R}$. Let $f$ be a single-valued self-mapping on $C$.

**Example 1.** Let $C = [0, \infty)$ and let $f$ be a single-valued self-mapping on $C$ satisfying

(i) $f(0) = 12$;
(ii) $0 \leq f(x) \leq x$, for $0 < x \leq 6$ and $f(0, 6] = [0, 6]$;
(iii) $f(x) \geq x$, for $6 \leq x < \infty$ and $f[6, \infty) = [6, \infty)$.

Then
(I) $f$ satisfies all conditions (a, b$_1$, c$_1$) in Theorem 1;
(II) $f$ has at least one fixed point, 6.

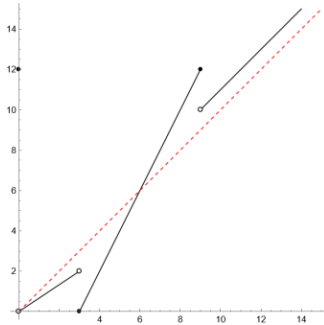

*Proof.* Conditions (i–iii) in this example show that $f$ satisfies condition (a) in Theorem 1. We next show that $f$ satisfies conditions (b$_1$). For any finite subset $\{x_1, x_2, \ldots, x_n\} \subseteq C$, let $u$ be an arbitrary convex combination of $x_1, x_2, \ldots, x_n$ with $u = \sum_{i=1}^{n} \alpha_i x_i$, for some positive numbers $\alpha_1, \alpha_2, \ldots, \alpha_n$ with $\sum_{i=1}^{n} \alpha_i = 1$. Suppose $0 \leq x_1 < x_2 < \ldots < x_n < \infty$. It implies that $0 \leq x_1 < u < x_n < \infty$. Then, the proof is divided to two cases:

Case 1. $0 \leq x_1 < u \leq 6$. By condition (ii), we have

$$0 < u - x_1 \leq u - f(x_1), \text{ for any } u \text{ with } 0 < x_1 < u \leq 6;$$

and
$$0 < u - 0 = u \leq 12 - u = f(0) - u, \text{ for any } u \text{ with } 0 = x_1 < u \leq 6.$$

Case 2. $6 < u < x_n < \infty$. By condition (iii), we have
$$0 < x_n - u \leq f(x_n) - u, \text{ for any } u \text{ with } 6 < u < x_n < \infty.$$

Then, we show that $f$ satisfies condition ($c_1$) in Theorem 1. Take any point $x^*$ with $x^* > 6$. Since $f(x^*) \geq x^*$, we have
$$\{y \in C : |x^* - y| \leq |f(x^*) - y|\} = \left[0, \tfrac{x^* + f(x^*)}{2}\right].$$
It is a compact subset in $C$. □

**Example 2.** Let $C = [0, \infty)$. Define two continuous functions φ, ψ on $C$ as follows:

$$\varphi(x) = \begin{cases} \frac{1}{2}x, & \text{for } x \in [0, 4]; \\ 3x - 10, & \text{for } x \in (4, 6]; \\ \frac{3}{2}x - 1, & \text{for } x \in (6, \infty), \end{cases}$$

and

$$\psi(x) = \begin{cases} \frac{3}{4}x, & \text{for } x \in [0, 4]; \\ 2x - 5, & \text{for } x \in (4, 6]; \\ \frac{5}{4}x - \frac{1}{2}, & \text{for } x \in (6, \infty). \end{cases}$$

Φ and ψ satisfy the following conditions:

(i) $0 < \varphi(x) < \psi(x) < x$, for $0 < x < 5$;
(ii) $\varphi(x) > \psi(x) > x$, for $5 < x < \infty$.

Then, we define $f: C \to C$ by
$$f(x) = \begin{cases} \varphi(x), & \text{for } 0 \leq x < \infty \text{ and } x \text{ is rational}; \\ \psi(x), & \text{for } 0 \leq x < \infty \text{ and } x \text{ is irrational}. \end{cases}$$
Then, we have

(I) $f$ satisfies all conditions (a, $b_1$, $c_1$) in Theorem 1;
(II) $f$ has two fixed points, 0 and 5;
(III) $f$ is discontinuous at every point in $(0, \infty)\setminus\{5\}$.

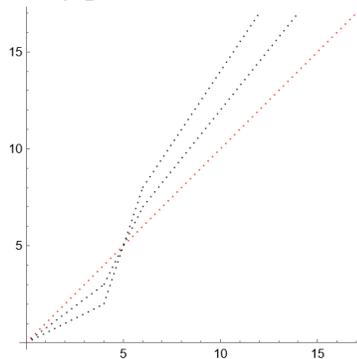

*Proof.* Proof of (a) in Theorem 1. Let $f(C)$ be the range of $f$. We have $\{0, 5\} \subseteq f(C)$. For any $b \in (0, \infty) \setminus \{5\}$, the horizontal line $y = b$ intersects with the curves $y = \varphi(x)$ and $y = \psi(x)$ at exactly one point, respectively. Suppose that $(x_1, b)$ and $(x_2, b)$ are the intersections of the line $y = b$ and the curves $y = \varphi(x)$ and $y = \psi(x)$, respectively. We have that if $b$ is rational, then both $x_1$ and $x_2$ are rational. It follows that $f(x_1) = \varphi(x_1) = b$. If $b$ is irrational, then both $x_1$ and $x_2$ are irrational. It follows that $f(x_2) = \psi(x_2) = b$. It proves (a). Similar to the proof of Example 1, we can show that $f$ satisfies condition ($b_1$) in Theorem 1.

Finally we show that $f$ satisfies condition ($c_1$) in Theorem 1. Take point $x^* = 6$. By $f(6) = 8$, we have
$$\{y \in C: |6 - y| \leq |f(6) - y|\} = \{y \in C: |6 - y| \leq |8 - y|\} = [0, 7].$$

It is a compact subset in $C$. □

### 2.3 Counter examples regarding to Theorem 1

Next, we give some counter examples to show that every condition in Theorem 1 is necessary for the considered mapping to have a fixed point.

**Example 3**. Let $C = [0, \infty)$ and let $f: C \to C$ be the function given in Example 2. Based on $f$, we define a function $g: C \to C$ by

$$g(x) = \begin{cases} f(x), & \text{for } x \neq 0, 5; \\ 10, & \text{for } x = 0; \\ 8, & \text{for } x = 5. \end{cases}$$

Then $g$ satisfies condition ($b_1$, $c_1$) but not (a) in Theorem 1 and $g$ has no fixed point.

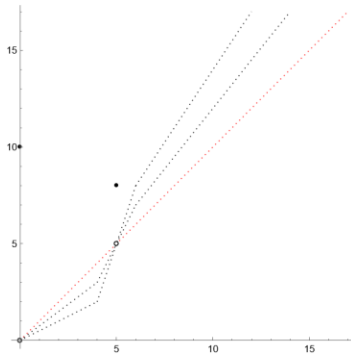

**Example 4**. Let $C = [0, 10]$. Let $f: C \to C$ be given in Example 5.

$$f(x) = \begin{cases} 0, & \text{for } x = 10; \\ \frac{4}{5}x, & \text{for } x \in (0, 5); \\ \frac{6}{5}x - 2, & \text{for } x \in [5, 10); \\ 10, & \text{for } x = 0. \end{cases}$$

Then *f* satisfies condition (a, c₁) but not (b₁) in Theorem 1 and *f* has no fixed point.

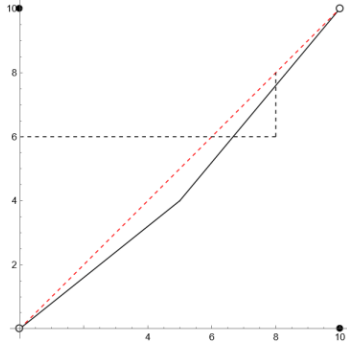

*Proof.* We only prove that *f* does not satisfy condition (b₁) in Theorem 1. Take $n = 2$ with $x_1 = 0$ and $x_2 > 5$. We can take $x_2$ very close to 10 such that there is *u* with $5 < u < x_2$ satisfying $u < f(x_2)$. Then we have
$$0 < u - x_1 = u > 10 - u = f(0) - u;$$
and
$$x_2 - u > f(x_2) - u > 0.$$

It shows that *f* does not satisfy condition (b₁) in Theorem 1. □

**Example 5.** Let $C = (-\infty, \infty)$ and define a linear function $f(x) = x - 1$. Then *f* satisfies conditions (a, b₁) but not (c₁) in Theorem 1 and *f* has no fixed point.

*Proof.* It is clear that *f* satisfies condition (a). Similar to the proof of Example 1, we can show that *f* satisfies condition (b₁) in Theorem 1. From the function $f(x) = x - 1$, we have

$$\{y \in C: |x - y| \leq |f(x) - y|\}$$
$$= \{y \in C: |x - y| \leq |x - 1 - y|\}$$
$$= \left[x - \frac{1}{2}, \infty\right), \text{ for any } x \in (-\infty, \infty).$$

It proves that *f* does not satisfy condition (c₁). □

### 3. The second fixed point theorem in this paper

The proof of the following theorem is very similar to the proof of Theorem 1. Hence, we only give a sketch proof of the following Theorem 3.

**Theorem 3.** *Let $(X, d)$ be a metric vector space and let C be a nonempty closed and convex subset of X. Let $f: C \to C$ be a single-valued mapping. Suppose that f satisfies the following conditions*:

(a) $f: C \to C$ *is onto*;

(b₂) *For any finite subset $\{x_1, x_2, \ldots, x_n\} \subseteq C$ and for arbitrary convex combination u of $x_1, x_2, \ldots, x_n$ with $u = \sum_{i=1}^{n} \alpha_i x_i$, in which $\alpha_1, \alpha_2, \ldots, \alpha_n$ are positive with $\sum_{i=1}^{n} \alpha_i = 1$,*

we have

$$\max\{d(f(x_j), u) - d(f(x_j), x_j): j = 1, 2, \ldots, n\} \geq 0;$$

(c$_2$) *There is $x^* \in C$ such that the following subset of $C$ is compact*

$$\{y \in C: d(f(x^*), x^*) \leq d(f(x^*), y)\},$$

*Then f has a fixed point.*

*Proof.* We define a set-valued mapping $G: C \to 2^C$ by

$$G(x) = \{y \in C: d(f(x), x) \leq d(f(x), y)\}, \text{ for every } x \in C.$$

Similar to the proof of Theorem 1, by condition (b$_2$), we can show that $G: C \to 2^C$ is a KKM mapping with nonempty closed values on $C$. By condition (c$_2$) and by using Fan-KKM Theorem, we obtain $\cap_{x \in C} G(x) \neq \emptyset$. Then, taking any $y_0 \in \cap_{x \in C} G(x)$, we have

$$d(f(x), x) \leq d(f(x), y_0), \text{ for every } x \in C. \tag{2}$$

By condition (a) in this theorem, for the arbitrarily selected $y_0 \in \cap_{x \in C} G(x) \subseteq C$ satisfying (2), there is $x_0 \in C$ such that $f(x_0) = y_0$. Substituting $x$ by $x_0$ in (2) gets

$$d(f(x_0), x_0) \leq d(f(x_0), y_0) = d(y_0, y_0) = 0.$$

It implies that $x_0$ is a fixed point of $f$, which proves this theorem. □

**Corollary 4.** *Let $(X, d)$ be a metric vector space and let $C$ be a nonempty compact and convex subset of $X$. Let $f: C \to C$ be a single-valued mapping. If $f$ satisfies conditions* (a) *and* (b$_2$) *in Theorem* 3, *then $f$ has a fixed point.*

Similar to the examples regarding to Theorem 1, we give three examples below for Corollary 4. The readers interested in these topics can construct more examples and counter examples to demonstrate Theorem 3 and Corollary 4.

**Example 6**. Let $C = [0, 10]$. Let $f: C \to C$ be defined with $f(C) = C$ and

$$f(x) \leq x, \text{ for all } x \in C.$$

Then, $f$ satisfies conditions (a, b$_2$) in Corollary 4 and $f$ has at least two fixed points, 0 and 10.

*Proof.* $f(C) = C$ is assumed, which proves (a). Next, we prove that $f$ satisfies condition (b$_2$). For any finite subset $\{x_1, x_2, \ldots, x_n\} \subseteq C$, let $u$ be an arbitrary convex combination of $x_1, x_2, \ldots, x_n$ with $u = \sum_{i=1}^{n} \alpha_i x_i$, for some positive numbers $\alpha_1, \alpha_2, \ldots, \alpha_n$ with $\sum_{i=1}^{n} \alpha_i = 1$. Suppose $0 \leq x_1 < u < x_n \leq 10$. Then, we have

$$0 \leq x_1 - f(x_1) < u - f(x_1), \text{ for any } u \text{ with } 0 \leq x_1 < u < 10.$$

□

Similar to Example 6, we have

**Example 7.** Let $C = [0, 10]$. Let $f: C \to C$ be defined with $f(C) = C$ and

$$f(x) \geq x, \text{ for all } x \in C.$$

Then, $f$ satisfies all conditions (a) and (b$_2$) in Corollary 4 and $f$ has at least two fixed points, 0 and 10.

*Proof.* From the proof of Theorem 6, we have

$$0 \leq f(x_n) - x_n < f(x_n) - u, \text{ for any } u \text{ with } 0 < u < x_n \leq 10. \qquad \square$$

**Example 8.** Let $C = [0, 10]$. Define two continuous functions $\varphi, \psi$ on $[0, 10]$ as follows:

$$\varphi(x) = \begin{cases} \frac{3}{5}x, & \text{for } x \in [0,5]; \\ \frac{7}{5}x - 4, & \text{for } x \in [5, 10]. \end{cases}$$

and

$$\psi(x) = \begin{cases} \frac{2}{5}x, & \text{for } x \in [0,5]; \\ \frac{8}{5}x - 6, & \text{for } x \in [5, 10]. \end{cases}$$

Based on the above given continuous functions $\varphi, \psi$, we define a function $f: C \to C$ by

$$f(x) = \begin{cases} \varphi(x), & \text{for } 0 \leq x \leq 10 \text{ and } x \text{ is rational}; \\ \psi(x), & \text{for } 0 < x < 10 \text{ and } x \text{ is irrational}. \end{cases}$$

Then, we have

(I)   $f$ satisfies all conditions (a, b$_2$) in Corollary 4;
(II)  $f$ has two fixed points, 0 and 10.
(III) $f$ is discontinuous at every point in $[0, 10]\setminus\{0, 10\}$.

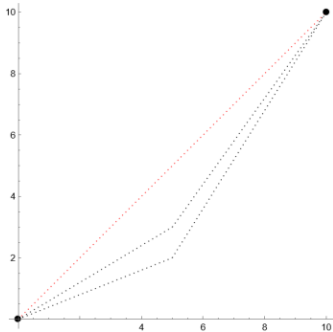

*Proof.* The proofs of (a) and (b$_2$) in this example are similar to the proofs of (a) in Example 2 and the proof of (b$_1$) in Example 6, respectively. $\qquad \square$

## 4. The third fixed point theorem in this paper

### 4.1 The third fixed point theorem in this paper

In this subsection, we use the Fan-KKM theorem to prove the third fixed point theorem in this paper. The way of the proof is different from the proofs of Theorems 1 and 3, which is said to be "a proof by indirectly applying the Fan-KKM theorem". For a self-mapping $f$ on a given set, let $\mathcal{F}(f)$ denote the collection of fixed point(s) of $f$.

**Theorem 5.** *Let $(X, d)$ be a metric vector space and let $C$ be a nonempty compact and convex subset of $X$. Let $f: C \to C$ be a single-valued mapping. Suppose that $f$ satisfies the following conditions*:

(a) $f: C \to C$ is onto;

(b$_3$) *For any finite subset $\{x_1, x_2, \ldots, x_n\} \subseteq C$ and for arbitrary convex combination $u$ of $x_1, x_2, \ldots, x_n$ with $u = \sum_{i=1}^{n} \alpha_i x_i$, in which $\alpha_1, \alpha_2, \ldots, \alpha_n$ are positive with $\sum_{i=1}^{n} \alpha_i = 1$, we have*
$$d(f(u), u) \leq \max\{d(f(x_j), u): j = 1, 2, \ldots, n\};$$

(c$_3$) $\{x \in C: d(f(x), x) \leq \beta\}$ *is closed, for any $\beta > 0$*;

*Then $\mathcal{F}(f)$ is a nonempty closed subset of $C$.*

*Proof.* $\mathcal{F}(f)$ is precisely defined by $\mathcal{F}(f) = \{x \in C: d(f(x), x) = 0\}$. Assume, by the way of contradiction, that
$$\mathcal{F}(f) = \emptyset. \tag{3}$$

Under the assumption (3), we show that there is $\delta > 0$ such that
$$\{x \in C: d(f(x), x) \leq \delta\} = \emptyset. \tag{4}$$

Assume by contradiction that (4) does not hold for any $\delta > 0$. That is,
$$\{x \in C: d(f(x), x) \leq \delta\} \neq \emptyset, \text{ for any } \delta > 0. \tag{5}$$

Let
$$E_m = \{x \in C: d(f(x), x) \leq \tfrac{1}{m}\}, \text{ for } m = 1, 2, \ldots.$$

By the assumption (5) and condition (c$_3$) in this theorem, $E_m$ is a nonempty closed subset of $C$, for $m = 1, 2, \ldots$. Then $\{E_m\}$ is a decreasing (with respect to the inclusion ordering) sequence of nonempty closed subsets of the given compact subset $C$ of $X$. It implies that
$$\cap_{m=1}^{\infty} E_m \neq \emptyset.$$
Since

$$\mathcal{F}(f) = \cap_{m=1}^{\infty} E_m, \tag{6}$$

it yields a contradiction to the assumption (3). Then (4) is proved. Hence, there is $\delta > 0$ that satisfies (4). Then we take $\frac{\delta}{2} > 0$ such that

$$d(f(x), x) \geq \frac{\delta}{2}, \text{ for every } x \in C. \tag{7}$$

With respect to this fixed $\frac{\delta}{2} > 0$ given in (7), we define a set-valued mapping $G: C \to 2^C$ by

$$G(x) = \{y \in C: d(f(x), y) \geq \frac{\delta}{2}\}, \text{ for every } x \in C.$$

By the assumption (7), we have $x \in G(x)$, for every $x \in C$. It follows that

$$C \supseteq G(x) \neq \emptyset, \text{ for every } x \in C.$$

Since, for every fixed $x \in C$, $d(f(x), \cdot)$ is a continuous function on $C$ (on $X$), it follows that, for every $x \in C$, $G(x)$ is a nonempty closed subset of $C$.

Next, we show that $G$ is a KKM mapping. To this end, for any finite subset $\{x_1, x_2, \ldots, x_n\} \subseteq C$, let $u$ be a linear convex combination of $x_1, x_2, \ldots, x_n$. We can suppose that there are positive numbers $\alpha_1, \alpha_2, \ldots, \alpha_n$ with $\sum_{i=1}^n \alpha_i = 1$ such that $u = \sum_{i=1}^n \alpha_i x_i$. By condition (b3) in this theorem, we have

$$d(f(u), u) = d(f(\sum_{i=1}^n \alpha_i x_i), \sum_{i=1}^n \alpha_i x_i) \leq \max\{d(f(x_j), \sum_{i=1}^n \alpha_i x_i): j = 1, 2, \ldots, n\}.$$

It implies that there must be an integer $k$ with $1 \leq k \leq n$ such that

$$d(f(u), u) = d(f(\sum_{i=1}^n \alpha_i x_i), \sum_{i=1}^n \alpha_i x_i) \leq d(f(x_k), \sum_{i=1}^n \alpha_i x_i).$$

By (7), it follows that

$$d(f(x_k), \sum_{i=1}^n \alpha_i x_i) \geq \frac{\delta}{2}.$$

That is,

$$u = \sum_{i=1}^n \alpha_i x_i \in G(x_k) \subseteq \cup_{1 \leq j \leq n} G(x_j).$$

It implies that $G: C \to 2^C$ is a KKM mapping with nonempty closed values in $C$. Since $C$ is compact, by Fan-KKM Theorem, we obtain $\cap_{x \in C} G(x) \neq \emptyset$. Taking any $y_0 \in \cap_{x \in C} G(x)$, we have

$$d(f(x), y_0) \geq \frac{\delta}{2}, \text{ for every } x \in C. \tag{8}$$

By condition (a) in this theorem, for the selected $y_0 \in C$ satisfying (8), there is $x_0 \in C$ such that $f(x_0) = y_0$. Substituting $x$ by $x_0$ in (8) gets

$$d(y_0, y_0) = d(f(x_0), y_0) \geq \frac{\delta}{2} > 0.$$

It is a contradiction. Hence (3) does not hold. That is, we must have $\mathcal{F}(f) \neq \emptyset$. By (6), $\mathcal{F}(f)$ is the intersection of a decreasing sequence of nonempty closed subsets of $C$. It follows that $\mathcal{F}(f)$ is closed. It proves this theorem. □

**Remarks.** The following (stronger) condition implies condition (b$_3$) in Theorem 5.

$$d(f(u), u) \leq \sum_{j=1}^{n} \alpha_j d(f(x_j), u).$$

### 4.2 Examples regarding to Theorem 5

In the following examples, we take $C = [0, 10]$. Let $f$ be a single-valued self-mapping on $C$.

**Example 9.**

$$f(x) = \begin{cases} -\frac{5}{3}x + 10, & \text{for } 0 \leq x < 3; \\ 1, & \text{for } x = 3; \\ 5, & \text{for } 3 < x < 7; \\ 9, & \text{for } x = 7; \\ -\frac{5}{3}x + \frac{50}{3}, & \text{for } 7 < x \leq 10. \end{cases}$$

Then we have

(i) $f$ satisfies all conditions in Theorem 5;

(ii) $f$ has one fixed point, 5;

(iii) $f$ is neither lower semi-continuous, nor upper semi-continuous.

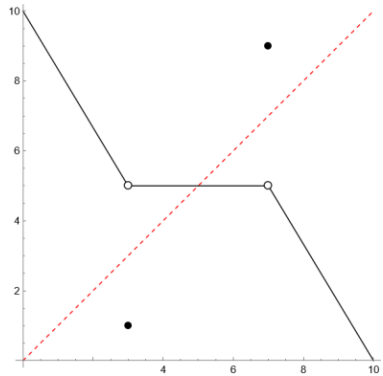

*Proof.* It is straight forward to check that $f$ satisfies conditions (a, c$_3$) in Theorem 5, and therefore it is omitted here. We only show that $f$ satisfies conditions (b$_3$). For any finite subset $\{x_1, x_2, \ldots, x_n\} \subseteq C$, let $u$ be an arbitrary convex combination of $x_1, \ldots, x_n$ with $u = \sum_{i=1}^{n} \alpha_i x_i$ for non-negative numbers $\alpha_1, \alpha_2, \ldots, \alpha_n$ with $\sum_{i=1}^{n} \alpha_i = 1$. Suppose $0 \leq x_1 < x_2 < \ldots < x_n \leq 10$. It is trivial for $u = x_1$ or $x_n$. So we assume that $0 \leq x_1 < u < x_n \leq 10$. Then, we have

$$|f(u) - u| \leq |f(x_1) - u|, \text{ for } x_1 < u \leq 5;$$

and

$$|f(u) - u| \leq |f(x_n) - u|, \text{ for } 5 < u < x_n. \qquad \square$$

Similar to Example 9, we construct more complicated examples. The proofs of the following two examples are very similar to the proof of Example 9, which are omitted here.

**Example 10.**

$$f(x) = \begin{cases} -\frac{5}{3}x + 10, & \text{for } 0 \leq x < 3; \\ 2x - 5, & \text{for } 3 \leq x \leq 7 \text{ and } x \text{ is rational}; \\ 5, & \text{for } 3 < x < 7 \text{ and } x \text{ is irrational}; \\ -\frac{5}{3}x + \frac{50}{3}, & \text{for } 7 < x \leq 10. \end{cases}$$

Then
  (i) $f$ satisfies all conditions in Theorem 5 and $f$ has one fixed point, 5.
  (ii) $f$ is continuous on $[0, 3) \cup (7, 10] \cup \{5\}$;
  (iii)   $f$ is discontinuous at every point in $[3, 7] \setminus \{5\}$;
  (iv)   $f$ is neither lower semi-continuous, nor upper semi-continuous on $C$.

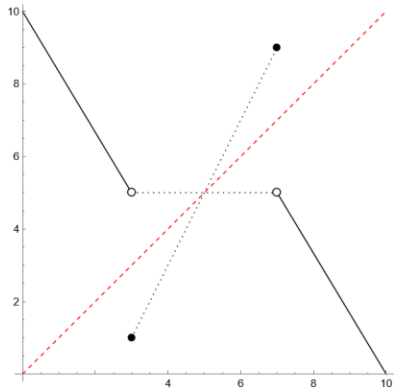

**Example 11.**

$$f(x) = \begin{cases} -\frac{5}{3}x + 10, & \text{for } 0 \leq x < 3; \\ 2x - 5, & \text{for } 3 \leq x \leq 7; \\ -\frac{5}{3}x + \frac{50}{3}, & \text{for } 7 < x \leq 10. \end{cases}$$

Then
  (i) $f$ satisfies all conditions in Theorem 5 and $f$ has one fixed point, 5;
  (ii) $f$ is continuous on $[0,10] \setminus \{3, 7\}$;
  (iii)   $f$ is neither lower semi-continuous, nor upper semi-continuous on $C$.

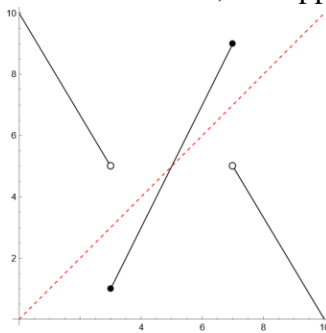

### 4.3 Counter examples regarding to Theorem 5

Now we provide some counterexamples to show that every condition in Theorem 5 is necessary for the considered mapping to have a fixed point.

**Example 12.**
$$f(x) = \begin{cases} 6, & \text{for } x \in [0, 5); \\ 4, & \text{for } x \in [5, 10). \end{cases}$$

Then $f$ satisfies all conditions (b$_3$) and (c$_3$) but not (a) in Theorem 5 and $f$ has no fixed point.

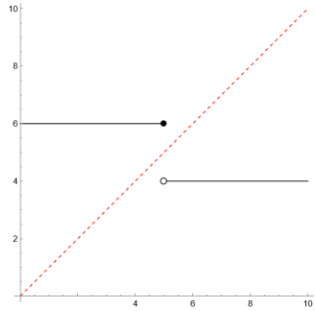

*Proof.* The proof is straight forward and it is omitted here. □

**Example 13.**
$$f(x) = \begin{cases} 0, & \text{for } x = 10; \\ \frac{4}{5}x, & \text{for } x \in (0, 5); \\ \frac{6}{5}x - 2, & \text{for } x \in [5, 10); \\ 10, & \text{for } x = 0. \end{cases}$$

Then $f$ satisfies all conditions (a) and (b$_3$) but not (c$_3$) in Theorem 5 and $f$ has no fixed point.

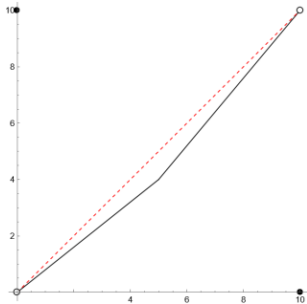

*Proof.* It is clear to see that $f$ satisfies condition (a). We next show that $f$ satisfies conditions (c). For any finite subset $\{x_1, x_2, \ldots, x_n\} \subseteq C$, let $u$ be an arbitrary convex combination of $x_1$, $x_2$, ..., $x_n$ with $u = \sum_{i=1}^{n} \alpha_i x_i$, for some positive numbers $\alpha_1, \alpha_2, \ldots, \alpha_n$ with $\sum_{i=1}^{n} \alpha_i = 1$. Suppose $0 \leq x_1 < x_2 < \ldots < x_n \leq 10$. It implies that $0 \leq x_1 < u < x_n \leq 10$. Then, we calculate

$$0 < u - f(u) < u - f(x_1), \text{ for any } 0 < x_1 < u < x_n \leq 10;$$

and
$$0 < u - f(u) < 10 - u, \text{ for any } 0 = x_1 < u < x_n \leq 10.$$

It shows that $f$ satisfies condition (b₃). Since $d(f(0), 0) = 10$, it implies that $\{x \in C: d(f(x), x) \leq \beta\}$ is not closed, for any $\beta$ with $0 < \beta < 1$. So, condition (c₃) is not satisfied. $f$ does not have a fixed point. □

**Example 14.**

$$f(x) = \begin{cases} -\frac{5}{3}x + 10, & \text{for } x \in [0, 3]; \\ x + 2, & \text{for } x \in (3, 5); \\ x - 2, & \text{for } x \in [5, 7]; \\ -\frac{5}{3}x + \frac{50}{3}, & \text{for } x \in (7, 10]. \end{cases}$$

Then $f$ satisfies conditions (a) and (c₃) but not (b₃) in Theorem 5 and $f$ has no fixed point.

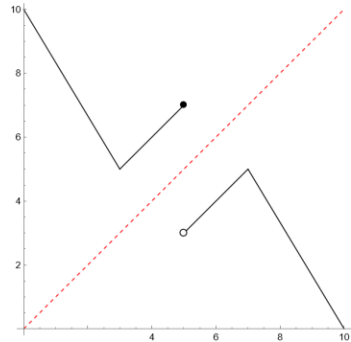

*Proof.* It is clear to see that $f$ satisfies conditions (a) and (c₃). Take $x_1 = 3$, $x_2 = 7$ and let $u = 5$. We have $f(x_1) - 5 = f(x_2) - 5 = 0$, and $f(u) - 5 = 2$. So, condition (b₃) is not satisfied. □

**Acknowledgments**

The author is very grateful to Professors Erdal Karapinar, Robert Mendris, Sehie Park, Adrian Petrusel, and Hongkun Xu for their insightful comments and suggestions in the development stage of this work.